\documentclass[10pt,reqno]{amsart}
\usepackage[T1]{fontenc}
\usepackage{textcomp}
\usepackage{mathtools}
\usepackage{mathptmx}
\usepackage{geometry}
\usepackage{enumitem}
\usepackage[expansion=false]{microtype}
\usepackage{xcolor}
\usepackage{hyperref}
\usepackage{amsmath}
\usepackage{amsfonts}
\usepackage{amssymb}
\usepackage{amsthm}
\usepackage{graphicx}
\hypersetup{
colorlinks=true,
linkcolor={blue!60!black},
citecolor={blue!60!black},
urlcolor={blue!80!black},
}

\theoremstyle{plain}
\newtheorem{theorem}{Theorem}[section]
\newtheorem{proposition}[theorem]{Proposition}
\newtheorem{lemma}[theorem]{Lemma}
\newtheorem{corollary}[theorem]{Corollary}
\theoremstyle{definition}

\newtheorem{example}[theorem]{Example}
\newtheorem{remark}[theorem]{Remark}
\theoremstyle{plain}
\newtheorem{assumption}{Assumption}
\makeatletter
\renewenvironment{proof}[1][\proofname]{\par
  \pushQED{\qed}%
  \normalfont \topsep6\p@\@plus6\p@\relax
  \trivlist
  \item[\hskip\labelsep\bfseries#1\@addpunct{.}]\ignorespaces
}{%
  \popQED\endtrivlist\@endpefalse
}
\makeatother
\DeclareMathOperator{\Image}{Im}
\DeclareMathOperator{\Ker}{ker}
\let\c\textcommabelow 
\begin{document}
\title[A constructive inverse for the Jordan--von Neumann operator]{A constructive
inverse for the Jordan--von Neumann operator on abelian groups}

\author{Dorian Popa}
\address{Department of Mathematics, Technical University of Cluj-Napoca,
         28 Memorandumului Street, 400114 Cluj-Napoca, Romania.\newline
         ORCID: \href{https://orcid.org/0000-0001-6197-3138}{0000-0001-6197-3138}}
\email{popa.dorian@math.utcluj.ro}
\thanks{$^{*}$Corresponding author: Dorian Popa
(\href{mailto:popa.dorian@math.utcluj.ro}{popa.dorian@math.utcluj.ro}).}

\author{Mircea Dan Rus}
\address{Department of Mathematics, Technical University of Cluj-Napoca,
         28 Memorandumului Street, 400114 Cluj-Napoca, Romania.\newline
         ORCID: \href{https://orcid.org/0000-0001-6255-0241}{0000-0001-6255-0241}}
\email{rus.mircea@math.utcluj.ro}

\date{1 July 2026; references and declarations updated 17 November 2026}

\begin{abstract}
For an abelian group $G$ and a real Banach space $X$, we study the
inhomogeneous Jordan--von Neumann equation $\mathsf{Q}[f] = g$, where
$\mathsf{Q}$ is the quadratic difference operator of the parallelogram law. We
characterize the functions $g$ in the image of $\mathsf{Q}$ and, under a
summability hypothesis, recover a solution explicitly by an iteration on the
diagonal followed by a propagation argument; neither regularity of $g$ nor
divisibility of $G$ is required. As a consequence we obtain a sharp stability
estimate of Forti--G\u{a}vru\c{t}a type for the quadratic functional equation.
\end{abstract}
\maketitle

\section{Introduction}

\label{sec:intro}

The Jordan--von Neumann functional equation
\begin{equation}
\label{eq:JvN}f(x + y) + f(x - y) = 2 f(x) + 2 f(y),
\end{equation}
is the algebraic skeleton of the parallelogram identity and, by the celebrated
theorem of Jordan and von Neumann~\cite{JordanVonNeumann1935}, characterizes
those norms on a real vector space that arise from an inner product. For
real-valued functions of a real variable its solution space admits an explicit
description: by Kurepa~\cite{Kurepa1956,Kurepa1964}, every $f\colon
\mathbb{R}\to\mathbb{R}$ satisfying~(\ref{eq:JvN}) is of the form $f(x) = B(x,
x)$ for a unique $\mathbb{Q}$-bilinear symmetric form $B\colon\mathbb{R}
\times\mathbb{R}\to\mathbb{R}$, and under any of several mild regularity
hypotheses (e.g., continuity at one point, Lebesgue measurability, or
boundedness on a set of positive measure) this collapses to $f(x) = c\,x^{2}$
for some $c\in\mathbb{R}$. Comprehensive accounts are given in the monographs
of Acz\'{e}l and Dhombres~\cite{AczelDhombres1989} and Kuczma~\cite{Kuczma2009}.

In the present paper we work in the natural generality of the equation. Let
$(G, +)$ be an abelian group and let $X$ be a real Banach space. For a function
$f\colon G\to X$ we consider the operator
\begin{equation}
\label{eq:Qdef}\mathsf{Q}[f](x, y) := f(x + y) + f(x - y) - 2 f(x) - 2
f(y),\qquad x, y\in G,
\end{equation}
which is $\mathbb{R}$-linear in $f$ and whose kernel $\Ker\mathsf{Q}$ is exactly
the solution space of~(\ref{eq:JvN}) on $G$. The solutions of~(\ref{eq:JvN}) are
the \emph{quadratic functions}: every such $f\colon G\to X$ has the form $f(x) =
B(x, x)$ for a unique symmetric biadditive $B\colon G\times G\to X$ (see, e.g.,
Stetk{\ae}r~\cite{Stetkaer2013} or Sz\'{e}kelyhidi~\cite{Szekelyhidi1991}); the
collapse to $f(x) = c\,x^{2}$ is a feature of the real line and has no analogue
on a general group.

These works describe the kernel of $\mathsf{Q}$, the quadratic functions
themselves; we are concerned instead with its image, through the
\emph{inhomogeneous} version of~(\ref{eq:JvN}): for a prescribed $g\colon
G\times G\to X$, we ask when there exists $f\colon G\to X$ with
\begin{equation}
\label{eq:inhom}\mathsf{Q}[f] = g,
\end{equation}
and how $f$ is then determined by $g$. The question is natural from two
perspectives. From the standpoint of the stability theory initiated by
Hyers~\cite{Hyers1941} and Ulam~\cite{Ulam1960}---and developed for the
quadratic equation by Skof~\cite{Skof1983}, Cholewa~\cite{Cholewa1984} and
Czerwik~\cite{Czerwik1992}, the passage to an arbitrary control function being
due to Forti~\cite{Forti1980,Forti1987} and G\u{a}vru\c{t}a~\cite{Gavruta1994};
see also~\cite{HyersIsacRassias1998,Jung2011}---the datum $g$ is a controlled perturbation of~(\ref{eq:JvN}), and
recovering $f$ from $g$ inverts the perturbation map. From an algebraic standpoint, the question parallels the
classical inverse problem for the linear Cauchy operator
\begin{equation}
\label{eq:Ldef}\mathsf{L}[f](x, y) := f(x + y) - f(x) - f(y),\qquad x, y\in G.
\end{equation}
For $\mathsf{L}$ it is a standard fact that $g$
belongs to the image of $\mathsf{L}$ if and only if $g$ satisfies both the
$2$-cocycle identity
\begin{equation}
\label{eq:Lcocycle}g(x, y) + g(x + y, z) = g(x, y + z) + g(y, z),\qquad x, y,
z\in G,
\end{equation}
and is symmetric. Symmetry is necessary because
$\mathsf{L}[f]$ is always symmetric, and it cannot be omitted: the cocycle
identity alone does not force it, as Erd\H{o}s~\cite{Erdos1959} showed on the
real line by producing asymmetric solutions of~(\ref{eq:Lcocycle}) outside
$\Image\mathsf{L}$ through a Hamel-basis construction. Conversely, once $g$ is
symmetric no regularity hypothesis is needed in the present setting: a real
vector space is divisible as an abelian group, and the resulting vanishing of
the relevant obstruction makes every symmetric solution
of~(\ref{eq:Lcocycle}) a coboundary; we return to this in the closing section.
On $\mathbb{R}$, Prunescu~\cite{Prunescu2007} gave for $g\in C^{1}$ an explicit
constructive counterpart. For the Sincov functional equation, another equation
of cocycle type, the inhomogeneous problem was treated recently by M\u{a}du\c{t}a, Popa,
and S\^{\i}ngeorzan~\cite{MadutaPopaSingeorzan2026}, who obtained solvability
conditions and the best Ulam constant.

The analogue of~(\ref{eq:Lcocycle}) for the quadratic operator $\mathsf{Q}$ is
the three-variable identity
\begin{equation}
\label{eq:cocycle}g(x + y, z) + g(x - y, z) + 2 g(x, y) = g(x, y + z) + g(x, y
- z) + 2 g(y, z),
\end{equation}
holding for all $x, y, z\in G$ whenever $g\in\Image\mathsf{Q}$. Its derivation
is parallel to that of~(\ref{eq:Lcocycle}) but rests on the parallelogram
structure of $\mathsf{Q}$ rather than on the associativity of addition (see
Section~\ref{sec:elementary}). By analogy with the linear case, we call
(\ref{eq:cocycle}) the \emph{cocycle identity} for $\mathsf{Q}$. It is
convenient to record the failure of~(\ref{eq:cocycle}) by means of the
$\mathbb{R}$-linear operator $\mathsf{D}\colon\mathcal{F}(G^{2}, X)\to
\mathcal{F}(G^{3}, X)$,
\begin{equation*}
\mathsf{D}[g](x, y, z) := g(x + y, z) + g(x - y, z) + 2 g(x, y)
- g(x, y + z) - g(x, y - z) - 2 g(y, z),
\end{equation*}
so that~(\ref{eq:cocycle}) reads $\mathsf{D}[g] = 0$; here $\mathcal{F}(G^{n},
X)$ denotes the space of all functions from $G^{n}$ to $X$.

A second condition will enter the characterization. We say that $g$ is
\emph{quasi-symmetric} if
\begin{equation}
\label{eq:qsym}g(x, y) - g(y, x) = g(0, 0) - g(0, x - y),\qquad x, y\in G.
\end{equation}
Equivalently, the antisymmetric part of $g$ is prescribed by its axial trace
$t\mapsto g(0, t)$; when this trace is constant, (\ref{eq:qsym}) reduces to
ordinary symmetry $g(x, y) = g(y, x)$. Quasi-symmetry plays here the role that
symmetry plays for $\mathsf{L}$: it is necessary (Section~\ref{sec:elementary})
and, like symmetry there, must be imposed alongside the cocycle identity rather
than deduced from it. On a general abelian group~(\ref{eq:cocycle}) does
\emph{not} force quasi-symmetry---we exhibit a counterexample on $\mathbb{Z}^{2}$
in Section~\ref{sec:cohomology}---whereas on $\mathbb{R}$ a continuity
hypothesis does force it; both points are taken up there.

This paper is a companion to~\cite{PaicuPopaRus2026}: there the inhomogeneous
equation~(\ref{eq:inhom}) is solved on the real line for a datum of class $C^{2}$
by a closed-form integral formula, with no restriction on the growth of $g$,
whereas here we treat the general setting $f\colon G\to X$ under a convergence
hypothesis of Forti--G\u{a}vru\c{t}a type. Both constructions are explicit, and the
two regimes overlap but neither contains the other: \cite{PaicuPopaRus2026}
trades growth for smoothness on $\mathbb{R}$, while the present construction
trades smoothness for a summability condition valid on any abelian group with
values in a Banach space, where completeness of $X$ is what turns the
construction into a convergent series.

The results are of three kinds. The two conditions just introduced are
\emph{necessary}: every $g\in\Image\mathsf{Q}$ satisfies the cocycle
identity~(\ref{eq:cocycle}) and is quasi-symmetric~(\ref{eq:qsym})
(Theorem~\ref{thm:necessity}, Section~\ref{sec:elementary}). Under a summability
hypothesis of Forti--G\u{a}vru\c{t}a type they are also \emph{sufficient}: a
solution of $\mathsf{Q}[f] = g$ is then produced explicitly, as a convergent
series in the diagonal values of $g$, unique up to a quadratic summand
(Theorem~\ref{thm:inverse}, Section~\ref{sec:formula}). The same construction
yields, as an application, a sharp Ulam stability estimate of
Forti--G\u{a}vru\c{t}a type for the quadratic equation (Theorem~\ref{thm:stab},
Section~\ref{sec:stability}).

Sections~\ref{sec:elementary}--\ref{sec:diagonal} are purely algebraic and use
only the real vector-space structure of $X$; completeness of $X$ enters in
Sections~\ref{sec:descent}--\ref{sec:propagation}. The plan is as follows.
Section~\ref{sec:elementary} records the elementary identities satisfied by
$\mathsf{Q}$ and uses them to prove the necessity of the cocycle identity and of
quasi-symmetry. Sections~\ref{sec:axes}--\ref{sec:diagonal} reduce the equation
to a one-variable recurrence on the diagonal by isolating the odd part and the
even reduction of the datum. Section~\ref{sec:descent} constructs the even part
by iterating that recurrence, Section~\ref{sec:propagation} propagates the
diagonal solution to all of $G\times G$, and Section~\ref{sec:formula} assembles
the explicit formula and settles uniqueness. Section~\ref{sec:examples} gives
examples, Section~\ref{sec:stability} proves the stability application, and
Section~\ref{sec:cohomology} closes with a cohomological interpretation.

\section{Elementary properties of $\mathsf{Q}$ and necessity}

\label{sec:elementary}

The operator $\mathsf{Q}$ defined by~(\ref{eq:Qdef}) is $\mathbb{R}$-linear. We
use throughout the canonical even/odd decomposition $f = f_{e} + f_{o}$, with
$f_{e}(x) = \tfrac{1}{2}\bigl(f(x) + f(-x)\bigr)$ and $f_{o}(x) =
\tfrac{1}{2}\bigl(f(x) - f(-x)\bigr)$; by linearity $\mathsf{Q}[f] =
\mathsf{Q}[f_{e}] + \mathsf{Q}[f_{o}]$. The next lemma collects five elementary
identities, all expressed through $f_{o}$ and the value $f(0)$.

\begin{lemma}
\label{lem:elementary}For all $f\in\mathcal{F}(G, X)$ and all $x, y\in G$,
\begin{align}
\mathsf{Q}[f](x, y) - \mathsf{Q}[f](x, -y)  &  = -4\, f_{o}(y),
\label{eq:el-signy}\\
\mathsf{Q}[f](x, y) - \mathsf{Q}[f](y, x)  &  = 2\, f_{o}(x - y),
\label{eq:el-swap}\\
\mathsf{Q}[f](x, y) - \mathsf{Q}[f](-x, -y)  &  = 2\, \mathsf{Q}[f_{o}](x, y),
\label{eq:el-bothsign}\\
\mathsf{Q}[f](x, 0)  &  = -2\, f(0), \label{eq:el-axis1}\\
\mathsf{Q}[f](0, y)  &  = -2\, f_{o}(y) - 2\, f(0). \label{eq:el-axis2}
\end{align}
\end{lemma}

\begin{proof}
All five follow by direct expansion of~(\ref{eq:Qdef}). For~(\ref{eq:el-signy}),
the substitution $y\mapsto-y$ exchanges $f(x + y)$ and $f(x - y)$ and turns
$-2 f(y)$ into $-2 f(-y)$, so the difference is $-2 f(y) + 2 f(-y) = -4
f_{o}(y)$. For~(\ref{eq:el-swap}), the swap $(x, y)\mapsto(y, x)$ fixes $f(x +
y)$, $-2 f(x)$ and $-2 f(y)$ and turns $f(x - y)$ into $f(y - x)$, so the
difference is $f(x - y) - f(y - x) = 2 f_{o}(x - y)$. For~(\ref{eq:el-bothsign}),
$(x, y)\mapsto(-x, -y)$ sends $f(x + y)$ to $f(-x - y)$, $f(x - y)$ to $f(-x +
y)$, $-2 f(x)$ to $-2 f(-x)$ and $-2 f(y)$ to $-2 f(-y)$; collecting the
differences gives $2 f_{o}(x + y) + 2 f_{o}(x - y) - 4 f_{o}(x) - 4 f_{o}(y) =
2\,\mathsf{Q}[f_{o}](x, y)$. Identity~(\ref{eq:el-axis1}) is $f(x) + f(x) - 2
f(x) - 2 f(0) = -2 f(0)$, and~(\ref{eq:el-axis2}) follows from $\mathsf{Q}
[f](0, y) = f(y) + f(-y) - 2 f(0) - 2 f(y)$ together with $f(-y) - f(y) = -2
f_{o}(y)$.
\end{proof}

\begin{corollary}
\label{cor:evenodd}Let $f\in\mathcal{F}(G, X)$.

\begin{enumerate}
[label=\textup{(\roman*)},leftmargin=*]

\item \label{it:even} If $f$ is even, then for all $x, y\in G$
\begin{gather*}
\mathsf{Q}[f](x, y) = \mathsf{Q}[f](-x, y) = \mathsf{Q}[f](x, -y) =
\mathsf{Q}[f](y, x),\\
\mathsf{Q}[f](x, 0) = \mathsf{Q}[f](0, y) = -2 f(0).
\end{gather*}

\item \label{it:odd} If $f$ is odd, then for all $x, y\in G$
\[
\mathsf{Q}[f](-x, -y) = -\mathsf{Q}[f](x, y),\qquad\mathsf{Q}[f](x, 0) =
0,\qquad\mathsf{Q}[f](0, y) = -2 f(y).
\]
\end{enumerate}
\end{corollary}

\begin{proof}
\ref{it:even} Here $f_{o} = 0$. By~(\ref{eq:el-signy}), $\mathsf{Q}[f]$ is even
in the second variable, and by~(\ref{eq:el-swap}) it is symmetric; together
these give the full chain of equalities, and the axial values follow
from~(\ref{eq:el-axis1})--(\ref{eq:el-axis2}).\smallskip

\ref{it:odd} Here $f_{e} = 0$, so $f_{o} = f$ and $f(0) = 0$. The antisymmetry
follows from~(\ref{eq:el-bothsign}), which reads $\mathsf{Q}[f](x, y) -
\mathsf{Q}[f](-x, -y) = 2\,\mathsf{Q}[f](x, y)$. The axial values follow
from~(\ref{eq:el-axis1})--(\ref{eq:el-axis2}) with $f(0) = 0$.
\end{proof}

These identities yield at once the necessity of the two conditions.

\begin{theorem}
[Necessity]\label{thm:necessity}If $g = \mathsf{Q}[f]$ for some $f\colon G\to
X$, then $g$ satisfies the cocycle identity~(\ref{eq:cocycle}) and is
quasi-symmetric, that is, satisfies~(\ref{eq:qsym}).
\end{theorem}

\begin{proof}
Set $g = \mathsf{Q}[f]$. For the cocycle identity, both sides
of~(\ref{eq:cocycle}) expand to
\begin{multline*}
S(x, y, z) := f(x + y + z) + f(x + y - z) + f(x - y + z) + f(x -
y - z)\\
- 4 f(x) - 4 f(y) - 4 f(z).
\end{multline*}
For the left-hand side, sum
\begin{align*}
g(x + y, z)  &  = f(x + y + z) + f(x + y - z) - 2 f(x + y) - 2 f(z),\\
g(x - y, z)  &  = f(x - y + z) + f(x - y - z) - 2 f(x - y) - 2 f(z),\\
2 g(x, y)  &  = 2 f(x + y) + 2 f(x - y) - 4 f(x) - 4 f(y);
\end{align*}
the $\pm2 f(x + y)$ and $\pm2 f(x - y)$ terms cancel, leaving $S(x, y, z)$, and
the analogous expansion of $g(x, y + z) + g(x, y - z) + 2 g(y, z)$ produces the
same $S(x, y, z)$. Hence $\mathsf{D}[g] = 0$. For quasi-symmetry,
(\ref{eq:el-swap}) gives $g(x, y) - g(y, x) = 2 f_{o}(x - y)$, while
(\ref{eq:el-axis2}) and~(\ref{eq:el-axis1}) give $g(0, t) = -2 f_{o}(t) - 2
f(0)$ and $g(0, 0) = -2 f(0)$, so $g(0, t) - g(0, 0) = -2 f_{o}(t)$. Taking
$t = x - y$ yields $2 f_{o}(x - y) = g(0, 0) - g(0, x - y)$,
which is~(\ref{eq:qsym}).
\end{proof}

In the special case $G = X = \mathbb{R}$, the cocycle identity and the elementary
identities of this section are due to~\cite{PaicuPopaRus2026}; the same arguments
apply verbatim on an arbitrary abelian group.

\begin{remark}
[Operatorial interpretation]\label{rem:operatorial}The cocycle
identities~(\ref{eq:Lcocycle}) and~(\ref{eq:cocycle}) are instances of a single
identity. For any linear defect operator $T\colon\mathcal{F}(G, X)\to\mathcal{F}
(G^{2}, X)$, applying $T$ to the two coordinate slices of $T[f]$ yields two
three-variable defects, and the \emph{iterated cocycle identity}
\begin{equation*}
T\bigl[\, T[f](\,\cdot\,, z)\,\bigr](x, y) = T\bigl[\,
T[f](x, \,\cdot\,)\,\bigr](y, z),\qquad x, y, z\in G,
\end{equation*}
asserts that they coincide; this is the cocycle condition for $T$, and gives the
standard $2$-cocycle identity for $T = \mathsf{L}$ and~(\ref{eq:cocycle}) for $T
= \mathsf{Q}$. This operatorial interpretation originates in~\cite{PaicuPopaRus2026}.
Casting the cocycle in this form---through the action of $T$ on the coordinate
slices of $g$---is what lets the linear and quadratic inverse problems be handled
in parallel; in Section~\ref{sec:cohomology} the same formulation yields a single
symmetry criterion specialising to the Erd\H{o}s--Acz\'{e}l theorem for $\mathsf{L}$
and to the quasi-symmetry of a continuous datum for $\mathsf{Q}$.
\end{remark}

\begin{remark}
[The quadratic kernel]\label{rem:kernel}In the introduction $\Ker\mathsf{Q}$ was
described through \emph{symmetric} biadditive maps; the symmetry is no
restriction. For any biadditive $B\colon G\times G\to X$ the antisymmetric part
vanishes on the diagonal, so $f(x) := B(x, x)$ equals $B_{\mathrm{s}}(x, x)$ with
$B_{\mathrm{s}}(x, y) = \tfrac{1}{2}\bigl(B(x, y) + B(y, x)\bigr)$ symmetric; in
particular $\mathsf{Q}[f] = 0$. A non-symmetric $B$ thus contributes nothing to
$\Ker\mathsf{Q}$ beyond its symmetrization.
\end{remark}

\section{Structural relations on the coordinate axes}

\label{sec:axes}

Throughout Sections~\ref{sec:axes}--\ref{sec:formula}, where $g$ is the datum to
be inverted, the following hypothesis is in force without being repeated.

\begin{assumption}
\label{ass:cocycle}$g\colon G\times G\to X$ is a fixed function satisfying the
cocycle identity~(\ref{eq:cocycle}).
\end{assumption}

\noindent Any further hypothesis on $g$ is stated explicitly where it is first
needed. We begin by extracting the relations satisfied by $g$ on the coordinate
axes, the foundation of the construction in the later sections.

\begin{lemma}
[Behaviour on the first axis]\label{lem:axis1}$g(x, 0) = g(0, 0)$ for all
$x\in G$; that is, $g$ is constant on the axis $G\times\{0\}$.
\end{lemma}

\begin{proof}
Setting $z = 0$ in~(\ref{eq:cocycle}), the right-hand side becomes $2 g(x, y) +
2 g(y, 0)$ and the left-hand side reads $g(x + y, 0) + g(x - y, 0) + 2 g(x, y)$.
Cancelling $2 g(x, y)$ gives
\begin{equation}
\label{eq:axis-aux}g(x + y, 0) + g(x - y, 0) = 2 g(y, 0).
\end{equation}
With $x = 0$ this yields $g(y, 0) + g(-y, 0) = 2 g(y, 0)$, so $g(-y, 0) = g(y,
0)$; that is, $y\mapsto g(y, 0)$ is even. Interchanging $x$ and $y$
in~(\ref{eq:axis-aux}) gives $g(x + y, 0) + g(y - x, 0) = 2 g(x, 0)$, whose
left-hand side equals that of~(\ref{eq:axis-aux}) by evenness. Hence $2 g(y, 0)
= 2 g(x, 0)$ for all $x, y$, so $g(\cdot, 0)$ is constant, with value $g(0, 0)$.
\end{proof}

\begin{lemma}
[Two key identities]\label{lem:axis2}For all $x, y, z\in G$,
\begin{align}
g(x, z) - g(x, -z)  &  = 2\bigl(g(0, z) - g(0, 0)\bigr), \label{eq:axis-signz}\\
g(-y, z) - g(y, z)  &  = g(0, y + z) + g(0, y - z) - 2 g(0, y).
\label{eq:axis-signy}
\end{align}
In particular, $z\mapsto g(0, z) - g(0, 0)$ is odd: $g(0, -z) + g(0, z) = 2
g(0, 0)$ for all $z\in G$.
\end{lemma}

\begin{proof}
For~(\ref{eq:axis-signz}), set $y = 0$ in~(\ref{eq:cocycle}): the left-hand side
is $g(x, z) + g(x, z) + 2 g(x, 0) = 2 g(x, z) + 2 g(0, 0)$ by
Lemma~\ref{lem:axis1}, while the right-hand side is $g(x, z) + g(x, -z) + 2
g(0, z)$. Equating and rearranging gives~(\ref{eq:axis-signz}). For
~(\ref{eq:axis-signy}), set $x = 0$ in~(\ref{eq:cocycle}): the left-hand side is
$g(y, z) + g(-y, z) + 2 g(0, y)$ and the right-hand side is $g(0, y + z) + g(0,
y - z) + 2 g(y, z)$; isolating $g(-y, z) - g(y, z)$ gives~(\ref{eq:axis-signy}).
Finally, setting $x = 0$ in~(\ref{eq:axis-signz}) yields $g(0, z) - g(0, -z) = 2
g(0, z) - 2 g(0, 0)$, i.e.\ $g(0, -z) - g(0, 0) = -\bigl(g(0, z) - g(0, 0)
\bigr)$, the asserted oddness.
\end{proof}

Identities~(\ref{eq:axis-signz}) and~(\ref{eq:axis-signy}) encode all the
asymmetry of $g$ forced by the cocycle identity; they are the main tools of the
next sections.

\section{Reduction to a one-variable equation}

\label{sec:diagonal}

We now begin the construction of a solution of~(\ref{eq:inhom}). By the linearity
of $\mathsf{Q}$ and the decomposition $f = f_{e} + f_{o}$, recovering $f$ amounts
to recovering $f_{e}$ and $f_{o}$ separately. The next proposition produces a
canonical candidate for the odd part, expressed through the values of $g$ on
the axis $\{0\}\times G$.

\begin{proposition}
\label{prop:odd}Define $\mu\colon G\to X$ by
\begin{equation}
\label{eq:phidef}\mu(x) := -\tfrac{1}{2}\bigl(g(0, x) - g(0, 0)\bigr).
\end{equation}
Then:

\begin{enumerate}
[label=\textup{(\roman*)},leftmargin=*]

\item \label{it:phi-odd} $\mu$ is odd;

\item \label{it:phi-Q} for all $x, y\in G$, $\mathsf{Q}[\mu]$ is given by
\begin{equation}
\label{eq:Qphi}\mathsf{Q}[\mu](x, y) = -\tfrac{1}{2}\bigl(g(0, x + y) +
g(0, x - y) - 2 g(0, x) - 2 g(0, y)\bigr) - g(0, 0),
\end{equation}
and in particular $\mathsf{Q}[\mu](x, 0) = 0$ and $\mathsf{Q}[\mu](0,
y) = g(0, y) - g(0, 0)$;

\item \label{it:phi-unique} if $f\colon G\to X$ satisfies $\mathsf{Q}[f] = g$,
then $f_{o} = \mu$.
\end{enumerate}
\end{proposition}

\begin{proof}
\ref{it:phi-odd} By Lemma~\ref{lem:axis2}, $t\mapsto g(0, t) - g(0, 0)$ is odd,
hence so is $\mu$.\smallskip

\ref{it:phi-Q} Writing $g(0, \cdot)$ for the function $t\mapsto g(0, t)$ and
$\mathbf{1}$ for the constant function $1$, linearity of $\mathsf{Q}$ gives
$\mathsf{Q}[\mu] = -\tfrac{1}{2}\mathsf{Q}[g(0, \cdot)] + \tfrac{1}{2} g(0,
0)\,\mathsf{Q}[\mathbf{1}]$; since $\mathsf{Q}[\mathbf{1}]\equiv-2$, expanding
$\mathsf{Q}[g(0, \cdot)]$ from~(\ref{eq:Qdef}) yields~(\ref{eq:Qphi}). The
axial values follow from Corollary~\ref{cor:evenodd}\,\ref{it:odd} applied to
the odd function $\mu$: $\mathsf{Q}[\mu](x, 0) = 0$ and $\mathsf{Q}
[\mu](0, y) = -2\mu(y) = g(0, y) - g(0, 0)$.\smallskip

\ref{it:phi-unique} Suppose $\mathsf{Q}[f] = g$. By~(\ref{eq:el-axis2}), $g(0,
y) = -2 f_{o}(y) - 2 f(0)$, and by~(\ref{eq:el-axis1}), $g(0, 0) = -2 f(0)$.
Subtracting gives $g(0, y) - g(0, 0) = -2 f_{o}(y)$, that is, $f_{o} = \mu$.
\end{proof}

Having isolated the odd part, we remove it: with $\mu$ as in
Proposition~\ref{prop:odd}, the \emph{even reduction} of $g$ is
\begin{equation}
\label{eq:gedef}g_{e} := g - \mathsf{Q}[\mu].
\end{equation}
Substituting~(\ref{eq:Qphi}) gives the closed form
\begin{equation}
\label{eq:geclosed}g_{e}(x, y) = g(x, y) + \tfrac{1}{2} g(0, x + y) +
\tfrac{1}{2} g(0, x - y) - g(0, x) - g(0, y) + g(0, 0).
\end{equation}

\begin{proposition}
\label{prop:residual}The even reduction $g_{e}$ satisfies the cocycle identity,
is even in each variable separately, and is constant on the coordinate axes:
\begin{align*}
g_{e}(x, y) &= g_{e}(-x, y) = g_{e}(x, -y),\\
g_{e}(x, 0) &= g_{e}(0, y) = g(0, 0),
\end{align*}
for all $x, y\in G$.
\end{proposition}

\begin{proof}
\emph{Cocycle.} By Theorem~\ref{thm:necessity}, $\mathsf{Q}[\mu]\in
\Ker\mathsf{D}$; since $g\in\Ker\mathsf{D}$ and $\mathsf{D}$ is linear, $g_{e} =
g - \mathsf{Q}[\mu]$ lies in $\Ker\mathsf{D}$ as well.\smallskip

\noindent\emph{Axial values.} By Proposition~\ref{prop:odd}\,\ref{it:phi-Q},
$\mathsf{Q}[\mu](x, 0) = 0$ and $\mathsf{Q}[\mu](0, y) = g(0, y) - g(0,
0)$. Combined with Lemma~\ref{lem:axis1} this gives $g_{e}(x, 0) = g(0, 0)$ and
\[
g_{e}(0, y) = g(0, y) - \bigl(g(0, y) - g(0, 0)\bigr) = g(0, 0).
\]
\smallskip

\noindent\emph{Evenness in the second variable.} Since $\mu$ is odd, (\ref{eq:el-signy})
and~(\ref{eq:phidef}) give
\[
\mathsf{Q}[\mu](x, y) - \mathsf{Q}[\mu](x, -y) = -4\mu(y) = 2\bigl(g(0, y) -
g(0, 0)\bigr),
\]
while~(\ref{eq:axis-signz}) gives $g(x, y) - g(x, -y) = 2\bigl(g(0, y) - g(0,
0)\bigr)$. Subtracting yields
\[
g_{e}(x, y) = g_{e}(x, -y).
\]
\smallskip

\noindent\emph{Evenness in the first variable.} Since $\mu$ is odd, (\ref{eq:phidef}) gives
\begin{align*}
\mathsf{Q}[\mu](-x, y) - \mathsf{Q}[\mu](x, y) & = -2\mu(x + y) - 2\mu(x - y) +
4\mu(x)\\
& = g(0, x + y) + g(0, x - y) - 2 g(0, x),
\end{align*}
while~(\ref{eq:axis-signy}) gives
\[
g(-x, y) - g(x, y) = g(0, x + y) + g(0, x - y) - 2 g(0, x).
\]
Subtracting yields $g_{e}(-x, y) = g_{e}(x, y)$.
\end{proof}

It remains to address the symmetry of $g_{e}$, which is precisely what the
quasi-symmetry of $g$ encodes.

\begin{proposition}
\label{prop:qsym-equiv}$g$ is quasi-symmetric if and only if $g_{e}$ is
symmetric.
\end{proposition}

\begin{proof}
From~(\ref{eq:geclosed}),
\[
g_{e}(x, y) - g_{e}(y, x) = \bigl(g(x, y) - g(y, x)\bigr) + \tfrac{1}{2}
\bigl(g(0, x - y) - g(0, y - x)\bigr).
\]
By the oddness of $t\mapsto g(0, t) - g(0, 0)$ (Lemma~\ref{lem:axis2}), $g(0, y
- x) = 2 g(0, 0) - g(0, x - y)$, so the right-hand side equals $\bigl(g(x, y) -
g(y, x)\bigr) + \bigl(g(0, x - y) - g(0, 0)\bigr)$. This vanishes for all $x, y$
precisely when $g(x, y) - g(y, x) = g(0, 0) - g(0, x - y)$, that is, when $g$ is
quasi-symmetric.
\end{proof}

We add a second standing hypothesis, in force through
Section~\ref{sec:formula}.

\begin{assumption}
\label{ass:qsym}$g$ is quasi-symmetric.
\end{assumption}

\noindent By Proposition~\ref{prop:qsym-equiv} the even reduction $g_{e}$ is then
symmetric. The reduction announced in the section title now takes the following
form.

\begin{lemma}
[Reduction to the diagonal]\label{lem:reduction}A function $f\colon G\to X$
solves $\mathsf{Q}[f] = g$ if and only if its odd part is $\mu$ and its even part
$f_{e}$ solves $\mathsf{Q}[f_{e}] = g_{e}$. In that case $f = f_{e} + \mu$, and
$f_{e}$ satisfies $f_{e}(0) = -\tfrac{1}{2} g(0, 0)$ together with the
one-variable recurrence
\begin{equation}
\label{eq:diag2}f_{e}(2 x) - 4 f_{e}(x) = g_{e}(x, x) + \tfrac{1}{2} g(0, 0),
\qquad x\in G.
\end{equation}
\end{lemma}

\begin{proof}
By the linearity of $\mathsf{Q}$ and $f = f_{e} + f_{o}$, the equation
$\mathsf{Q}[f] = g$ splits into its odd and even parts:
Proposition~\ref{prop:odd}\,\ref{it:phi-unique} forces $f_{o} = \mu$, and then
$\mathsf{Q}[f_{e}] = g - \mathsf{Q}[\mu] = g_{e}$. For the recurrence, evaluating
$\mathsf{Q}[f_{e}] = g_{e}$ at $y = 0$ gives $-2 f_{e}(0) = g_{e}(x, 0) = g(0,
0)$ by~(\ref{eq:el-axis1}) and Proposition~\ref{prop:residual}, so $f_{e}(0) =
-\tfrac{1}{2} g(0, 0)$; at $y = x$, using $\mathsf{Q}[f_{e}](x, x) = f_{e}(2 x) +
f_{e}(0) - 4 f_{e}(x)$, it gives~(\ref{eq:diag2}).
\end{proof}

The problem is thus reduced to solving the single recurrence~(\ref{eq:diag2})
for an even unknown function. The next section does this in the Banach space $X$,
and the one after propagates the diagonal solution to all of $G\times G$.

\section{Construction of the even part}

\label{sec:descent}

Throughout this section and the next we impose, in addition to the standing
hypotheses, one further assumption on $g$.

\begin{assumption}
[Two-variable Forti--G\u{a}vru\c{t}a condition]\label{ass:gavruta}For all $x, y\in G$,
\begin{equation}
\label{eq:conv}\sum_{n\geq1}\frac{\bigl\|g(2^{n-1}x, 2^{n-1}y)\bigr\|}{4^{n}} <
\infty.
\end{equation}
\end{assumption}

\begin{remark}
\label{rem:gavruta-role}Condition~(\ref{eq:conv}) controls the growth of $g$
along the doubling orbit and is what makes the defining series of the
construction converge. Its force depends on the setting: for bounded data---in
particular on any finite or compact group---it holds with no further
hypothesis, whereas on $\mathbb{R}$ it confines the growth of $g$ to below the
quadratic rate, failing exactly at the critical exponent. The examples of
Section~\ref{sec:examples} display both sides, from recovery throughout the
subcritical range to the breakdown at the parallelogram-law defect.
\end{remark}

\noindent We construct an even function $\nu$ solving the diagonal
recurrence~(\ref{eq:diag2}), using the completeness of $X$.

Putting $y = x$ in~(\ref{eq:geclosed}) gives the diagonal of $g_{e}$,
\begin{equation}
\label{eq:gediag}g_{e}(x, x) = g(x, x) + \tfrac{1}{2} g(0, 2 x) - 2 g(0, x) +
\tfrac{3}{2} g(0, 0).
\end{equation}
The next lemma records the series that will both define and bound $\nu$.

\begin{lemma}
\label{lem:psi}For every $x\in G$ the series
\begin{equation}
\label{eq:psidef}\psi(x) := \sum_{n\geq1} 4^{-n}\bigl\|g_{e}(2^{n-1}x,
2^{n-1}x)\bigr\|
\end{equation}
converges, and
\begin{equation}
\label{eq:decay}\lim\limits_{n\to\infty} 4^{-n}\psi(2^{n}x) = 0.
\end{equation}
\end{lemma}

\begin{proof}
Each term on the right of~(\ref{eq:gediag}), evaluated at $2^{n-1}x$, is summable
against $4^{-n}$ by the Forti--G\u{a}vru\c{t}a condition~(\ref{eq:conv}) (its diagonal $y = x$
and its axis $x = 0$); hence $\psi(x) < \infty$. For~(\ref{eq:decay}), expanding
$\psi(2^{n}x)$ and reindexing $i = m + n$ give
\[
4^{-n}\psi(2^{n}x) = \sum_{m\geq1} 4^{-(m + n)}\bigl\|g_{e}(2^{m + n - 1}x, 2^{m
+ n - 1}x)\bigr\| = \sum_{i\geq n + 1} 4^{-i}\bigl\|g_{e}(2^{i - 1}x, 2^{i -
1}x)\bigr\|,
\]
the tail of the convergent series~(\ref{eq:psidef}), which tends to $0$.
\end{proof}

\begin{proposition}
\label{prop:descent}The function $\nu\colon G\to X$ defined by
\begin{equation}
\label{eq:fedef}\nu(x) := -\sum_{n\geq1} 4^{-n} g_{e}(2^{n-1}x, 2^{n-1}x)
- \tfrac{1}{6} g(0, 0)
\end{equation}
is even and satisfies
\begin{enumerate}
[label=\textup{(\roman*)},leftmargin=*]
\item $\nu(0) = -\tfrac{1}{2} g(0, 0)$;

\item the recurrence~(\ref{eq:diag2});

\item $\mathsf{Q}[\nu](x, x) = g_{e}(x, x)$ for all $x\in G$;

\item $\|\nu(x)\|\leq\psi(x) + \tfrac{1}{6}\|g(0, 0)\|$ for all $x$.
\end{enumerate}
\end{proposition}

\begin{proof}
\emph{Well-definedness and bound.} Since $\sum_{n\geq1} 4^{-n}\bigl\|g_{e}(2^{n-1}x,
2^{n-1}x)\bigr\| = \psi(x) < \infty$ by Lemma~\ref{lem:psi}, the series
in~(\ref{eq:fedef}) converges in $X$ by completeness, so $\nu$ is well defined;
the triangle inequality and $\sum_{n\geq1} 4^{-n} = \tfrac{1}{3}$ give
$\|\nu(x)\|\leq\psi(x) + \tfrac{1}{6}\|g(0, 0)\|$.\smallskip

\noindent\emph{Evenness.} As $g_{e}$ is even in each variable
(Proposition~\ref{prop:residual}), each summand in~(\ref{eq:fedef}) is unchanged
under $x\mapsto-x$, so $\nu(-x) = \nu(x)$.\smallskip

\noindent\emph{Value at $0$ and recurrence.} With $d(x) := g_{e}(x, x) + \tfrac{1}{2} g(0,
0)$ and $\sum_{n\geq1} 4^{-n} = \tfrac{1}{3}$, formula~(\ref{eq:fedef}) reads
$\nu(x) = -\sum_{n\geq1} 4^{-n} d(2^{n-1}x)$. Since $g_{e}(0, 0) = g(0, 0)$, we
have $d(0) = \tfrac{3}{2} g(0, 0)$ and hence $\nu(0) = -\tfrac{1}{3} d(0) =
-\tfrac{1}{2} g(0, 0)$. For the recurrence,
\[
4 \nu(x) = -\sum_{n\geq1} 4^{-(n-1)} d(2^{n-1}x) = -d(x) - \sum_{m\geq1} 4^{-m}
d(2^{m}x) = -d(x) + \nu(2 x),
\]
by reindexing $m = n - 1$; hence $\nu(2 x) - 4 \nu(x) = d(x)$, which
is~(\ref{eq:diag2}), and $\mathsf{Q}[\nu](x, x) = \nu(2 x) + \nu(0) - 4
\nu(x) = d(x) - \tfrac{1}{2} g(0, 0) = g_{e}(x, x)$.
\end{proof}

\section{Propagation to the full equation}

\label{sec:propagation}

Let $\mu$ be as in Proposition~\ref{prop:odd},
$\nu$ as in Proposition~\ref{prop:descent}, and set
\begin{equation}
\label{eq:fdef}f := \nu + \mu.
\end{equation}
Since $g_{e} = g - \mathsf{Q}[\mu]$, we have $\mathsf{Q}[f] = \mathsf{Q}
[\nu] + \mathsf{Q}[\mu] = g$ as soon as $\mathsf{Q}[\nu] = g_{e}$. By
Proposition~\ref{prop:descent} the two sides already agree on the diagonal; it
remains to extend the agreement to all of $G\times G$. Put
\begin{equation*}
E := \mathsf{Q}[\nu] - g_{e}.
\end{equation*}

\begin{lemma}
\label{lem:Eprops}The function $E$ lies in $\Ker\mathsf{D}$, is symmetric and
even in each variable, and vanishes on the diagonal and on the coordinate axes.
\end{lemma}

\begin{proof}
Both $\mathsf{Q}[\nu]$ and $g_{e}$ lie in $\Ker\mathsf{D}$---the former by
Theorem~\ref{thm:necessity}, the latter by Proposition~\ref{prop:residual}---so
$E\in\Ker\mathsf{D}$ by linearity. As $\nu$ is even,
Corollary~\ref{cor:evenodd}\,\ref{it:even} makes $\mathsf{Q}[\nu]$ symmetric
and even in each variable; $g_{e}$ has the same properties
(Proposition~\ref{prop:residual} and, for symmetry,
Proposition~\ref{prop:qsym-equiv}), so $E$ does too. On the diagonal $E(x, x) =
\mathsf{Q}[\nu](x, x) - g_{e}(x, x) = 0$ by Proposition~\ref{prop:descent}.
Finally $\mathsf{Q}[\nu](x, 0) = -2 \nu(0) = g(0, 0)$ by~(\ref{eq:el-axis1})
and $\mathsf{Q}[\nu](0, y) = -2 \nu(0) = g(0, 0)$ by~(\ref{eq:el-axis2})
($\nu$ even), while $g_{e} = g(0, 0)$ on the axes
(Proposition~\ref{prop:residual}); hence $E$ vanishes on the axes.
\end{proof}

\begin{lemma}
[Doubling relation]\label{lem:doubling}For every $x, y\in G$,
\begin{equation}
\label{eq:star}E(2 x, 2 y) = 2 E(x + y, x - y) + 8 E(x, y).
\end{equation}
\end{lemma}

\begin{proof}
We use the cocycle identity for $E$ together with $E(x, x) = 0$ and the
vanishing of $E$ on the axes. Taking $a = b = x$ in the cocycle (in the
form~(\ref{eq:cocycle}) with $g$ replaced by $E$) and using $E(0, z) = 0$ gives
\begin{equation}
\label{eq:doubA}E(2 x, z) = E(x, x + z) + E(x, x - z) + 2 E(x, z).
\end{equation}
Setting $z = 2 y$ in~(\ref{eq:doubA}),
\begin{equation}
\label{eq:doubA1}E(2 x, 2 y) = E(x, x + 2 y) + E(x, x - 2 y) + 2 E(x, 2 y).
\end{equation}
Next, taking $(a, b, c) = (x, y, y)$ in the cocycle and using $E(x, 0) = 0 =
E(y, y)$,
\begin{equation}
\label{eq:doubB}E(x, 2 y) = E(x + y, y) + E(x - y, y) + 2 E(x, y).
\end{equation}
Taking $(a, b, c) = (x, y, x + y)$ and using $E(x + y, x + y) = 0$, $E(x, -x) =
E(x, x) = 0$ and symmetry,
\begin{equation}
\label{eq:doubC1}E(x, x + 2 y) = E(x - y, x + y) + 2 E(x, y) - 2 E(x + y, y);
\end{equation}
replacing $y$ by $-y$ and using evenness gives
\begin{equation}
\label{eq:doubC2}E(x, x - 2 y) = E(x + y, x - y) + 2 E(x, y) - 2 E(x - y, y).
\end{equation}
Substituting~(\ref{eq:doubB}), (\ref{eq:doubC1}) and~(\ref{eq:doubC2})
into~(\ref{eq:doubA1}), the terms $\mp2 E(x + y, y)$ and $\mp2 E(x - y, y)$
cancel and the two terms $E(x - y, x + y) = E(x + y, x - y)$ combine,
leaving~(\ref{eq:star}).
\end{proof}

\begin{theorem}
\label{thm:propagation}Under Assumptions 1--3, $\mathsf{Q}[\nu] =
g_{e}$; consequently the function $f = \nu + \mu$
defined by~(\ref{eq:fdef}) satisfies $\mathsf{Q}[f] = g$.
\end{theorem}

\begin{proof}
Fix $x, y\in G$ and write $R(u, v) := (u + v, u - v)$, so that $R^{2}(u, v) =
(2 u, 2 v)$. Applying Lemma~\ref{lem:doubling} at the point $R^{k}(x, y)$, the
sequence $a_{k} := E\bigl(R^{k}(x, y)\bigr)$, $k\geq0$, in $X$ satisfies the
linear recurrence
\[
a_{k + 2} = 2 a_{k + 1} + 8 a_{k},\qquad k\geq0,
\]
whose characteristic roots are $4$ and $-2$; hence $a_{k} = 4^{k}\alpha +
(-2)^{k}\beta$ for some $\alpha, \beta\in X$ depending on $(x, y)$. Since
$R^{2 n}(x, y) = (2^{n}x, 2^{n}y)$,
\begin{equation}
\label{eq:Egrowth}E(2^{n}x, 2^{n}y) = a_{2 n} = 16^{n}\alpha + 4^{n}\beta,
\qquad n\geq0.
\end{equation}

We now show that $4^{-n}\bigl\|E(2^{n}x, 2^{n}y)\bigr\|\to0$. By the norm
bound of Proposition~\ref{prop:descent}, applied to
$w\in\{x + y, x - y, x, y\}$,
\begin{multline*}
4^{-n}\bigl\|\mathsf{Q}[\nu](2^{n}x, 2^{n}y)\bigr\|\\
\leq 4^{-n}\bigl(\psi(2^{n}(x + y)) + \psi(2^{n}(x - y)) + 2\psi(2^{n}x) +
2\psi(2^{n}y) + \|g(0, 0)\|\bigr),
\end{multline*}
which tends to $0$ by~(\ref{eq:decay}). On the other hand,~(\ref{eq:geclosed})
bounds $\|g_{e}(2^{n}x, 2^{n}y)\|$ by a fixed combination of $\|g(2^{n}x,
2^{n}y)\|$ and the axial norms $\|g(0, 2^{n}w)\|$, $w\in\{x + y, x - y, x, y\}$.
The convergence of~(\ref{eq:conv}) forces its general term to vanish,
\[
4^{-n}\bigl\|g(2^{n-1}u, 2^{n-1}v)\bigr\|\to0,\qquad\text{hence}\qquad
4^{-n}\bigl\|g(2^{n}u, 2^{n}v)\bigr\|\to0,
\]
for all $u, v\in G$ (the axes being the case $u = 0$), whence
$4^{-n}\|g_{e}(2^{n}x, 2^{n}y)\|\to0$. Together with the bound on
$\mathsf{Q}[\nu]$ above, this gives $4^{-n}\|E(2^{n}x, 2^{n}y)\|\to0$.

By~(\ref{eq:Egrowth}), $4^{-n}E(2^{n}x, 2^{n}y) = 4^{n}\alpha + \beta$; as its
norm tends to $0$, letting $n\to\infty$ forces $\alpha = 0$, and then $\beta = 0$.
Hence $E(x, y) = a_{0} = \alpha + \beta = 0$. As $(x, y)$ was arbitrary,
$E\equiv0$, i.e.\ $\mathsf{Q}[\nu] = g_{e}$.
\end{proof}

\section{The inverse: explicit formula and uniqueness}

\label{sec:formula}

Assembling Propositions~\ref{prop:odd} and~\ref{prop:descent} through
Theorem~\ref{thm:propagation} yields the constructive half of the
characterization, with an explicit formula for the solution.

\begin{theorem}
\label{thm:inverse}Let $g\colon G\times G\to X$ be a quasi-symmetric cocycle
(that is, satisfy~(\ref{eq:cocycle}) and~(\ref{eq:qsym})) obeying the two-variable
Forti--G\u{a}vru\c{t}a condition~(\ref{eq:conv}). Then the function $f\colon G\to X$ given
by
\begin{equation}
\label{eq:explicit}f(x) = -\sum_{n\geq1} 4^{-n}\, g(2^{n-1}x, 2^{n-1}x)
- \tfrac{1}{6} g(0, 0),\qquad x\in G,
\end{equation}
satisfies $\mathsf{Q}[f] = g$, and the solution set of $\mathsf{Q}[f] = g$ is $f
+ \Ker\mathsf{Q}$.
\end{theorem}

\begin{proof}
By Theorem~\ref{thm:propagation}, $f = \nu + \mu$ solves $\mathsf{Q}[f] = g$,
with $\mu$ from~(\ref{eq:phidef}) and $\nu$ the series~(\ref{eq:fedef}). Inserting
the diagonal value~(\ref{eq:gediag}) into~(\ref{eq:fedef}) and using
$\sum_{n\geq1} 4^{-n} = \tfrac{1}{3}$ gives, with $b_{n} := 4^{-n} g(0, 2^{n}x)$,
\[
\nu(x) = -\sum_{n\geq1} 4^{-n} g(2^{n-1}x, 2^{n-1}x) - \tfrac{1}{2}
\sum_{n\geq1}(b_{n} - b_{n-1}) - \tfrac{2}{3} g(0, 0).
\]
The middle sum telescopes to $-g(0, x)$, as $b_{0} = g(0, x)$ and $b_{n} = 4^{-n}
g(0, 2^{n}x)\to0$ by~(\ref{eq:conv}); hence $\nu(x) = -\sum_{n\geq1} 4^{-n}
g(2^{n-1}x, 2^{n-1}x) + \tfrac{1}{2} g(0, x) - \tfrac{2}{3} g(0, 0)$, and adding
$\mu$ cancels the axial term to give~(\ref{eq:explicit}). If $f_{1}$ also
solves $\mathsf{Q}[f_{1}] = g$, then $\mathsf{Q}[f_{1} - f] = 0$, so
$f_{1} - f\in\Ker\mathsf{Q}$ is quadratic; conversely every such $f_{1}$
solves the equation.
\end{proof}

Together with Theorem~\ref{thm:necessity}, this settles the inverse problem: for
$g$ satisfying the summability~(\ref{eq:conv}), membership in $\Image\mathsf{Q}$
is equivalent to the cocycle identity together with quasi-symmetry.

\begin{remark}
\label{rem:diagdom}Although~(\ref{eq:conv}) is a condition on all pairs $(x,
y)$, it follows from its diagonal alone when $g$ is dominated by its diagonal,
that is, when
\[
\|g(x, y)\|\leq\tfrac{1}{2}\bigl(\|g(x, x)\| + \|g(y, y)\|\bigr),\qquad x, y\in G.
\]
Then convergence of $\sum_{n\geq1} 4^{-n}\|g(2^{n-1}x, 2^{n-1}x)\|$ for all $x$
already gives~(\ref{eq:conv}). Such a bound
cannot be assumed for a prescribed datum, which is why~(\ref{eq:conv}) is stated
for all pairs. In the stability application of Section~\ref{sec:stability}, by
contrast, the control function is chosen by the user, and the standard choices
do satisfy the corresponding bound (Corollary~\ref{cor:diagdom}).
\end{remark}

\section{Examples}

\label{sec:examples}

We illustrate the construction from several sides: how its odd and even parts
combine; the real line, where it is instructive to compare with the smooth
regime of~\cite{PaicuPopaRus2026}; the parallelogram-law defect of a
non-Hilbertian norm, which displays at once the reach of the method and its
sharp limit; discrete, finite, and compact groups, where no smooth regime exists
and the membership condition turns purely algebraic; and an infinite-dimensional
target, where completeness of the target is essential. In each example we
exhibit a single solution $f$; by the linearity of $\mathsf{Q}$ the full solution
set is then $f + \Ker\mathsf{Q}$, so $f$ is determined uniquely modulo
$\Ker\mathsf{Q}$.

\begin{example}
[Constant datum]\label{ex:const}Let $g\equiv c$ be constant, $c\in X$. Then $g$
satisfies~(\ref{eq:cocycle}) and~(\ref{eq:qsym}) trivially,
and~(\ref{eq:conv}) holds. Here $g_{e}\equiv c$ and $\mu\equiv0$,
so~(\ref{eq:explicit}) gives $f\equiv-\tfrac{1}{2} c$, and indeed $\mathsf{Q}
[-\tfrac{1}{2} c] = c$.
\end{example}

\begin{example}
[An additive datum]\label{ex:additive}Let $a\colon G\to X$ be additive and put
$g(x, y) := a(y)$, independent of the first variable. By additivity
$\mathsf{Q}[a](x, y) = -2 a(y)$, so $g = \mathsf{Q}[-\tfrac{1}{2} a]$ lies in
$\Image\mathsf{Q}$ and is therefore a quasi-symmetric cocycle.
Condition~(\ref{eq:conv}) holds, as $\sum_{n\geq1} 4^{-n}\|a(2^{n-1}y)\| =
\tfrac{1}{2}\|a(y)\|$. Here $g(0, x) = a(x)$ gives $\mu(x) = -\tfrac{1}{2} a(x)$,
while $g_{e}\equiv0$ and hence $\nu\equiv0$; formula~(\ref{eq:explicit}) returns
$f = -\tfrac{1}{2} a$. Since an additive function is automatically odd, $a(-x) =
-a(x)$, this is a purely odd solution, with vanishing even part.
\end{example}

\begin{example}
[A non-smooth datum]\label{ex:abs}Let $G = \mathbb{R}$, $X = \mathbb{R}$, and
$g = \mathsf{Q}[h]$ with $h(x) = |x|$, so that $g(x, y) = |x + y| + |x - y| - 2
|x| - 2 |y|$. Then $g(0, x) = 0$ and $g_{e}(t, t) = g(t, t) = -2 |t|$, while the
linear growth of $g$ makes~(\ref{eq:conv}) hold. Formula~(\ref{eq:explicit})
returns
\[
f(x) = -\sum_{n\geq1} 4^{-n}\bigl(-2\,|2^{n-1}x|\bigr) = |x|\sum_{n\geq1}
2^{-n} = |x|,
\]
recovering $h$ exactly. As $h$ is continuous but not differentiable, this datum
lies outside the $C^{2}$ regime of~\cite{PaicuPopaRus2026}, illustrating the
complementarity of the two approaches.
\end{example}

\begin{example}
[A datum with both parts]\label{ex:mixed}On $G = X = \mathbb{R}$ take $h(x) =
|x| + \sin x$ and $g = \mathsf{Q}[h]$, whose even and odd parts are $|x|$ and
$\sin x$. From $g(0, x) = h(x) + h(-x) - 2 h(0) - 2 h(x) = -2\sin x$,
formula~(\ref{eq:phidef}) gives $\mu(x) = \sin x$, the odd part exactly; the
even reduction has diagonal $g_{e}(t, t) = -2 |t|$ by~(\ref{eq:gediag})---the
odd part being carried entirely by $\mu$---so $\nu(x) = |x|$ as in
Example~\ref{ex:abs}, and~(\ref{eq:explicit}) returns $f = |x| + \sin x = h$.
Here $g$ is genuinely non-symmetric, $g(x, y) - g(y, x) = 2\sin(x - y)$, so it
is quasi-symmetry, not symmetry, that is verified---a feature with no analogue
for the additive Cauchy operator.
\end{example}

\begin{example}
[A bounded datum]\label{ex:cos}For $h(x) = \cos x$ on $\mathbb{R}$ one has $g =
\mathsf{Q}[h]$ bounded, so~(\ref{eq:conv}) holds. Here
$g(0, x) = -2$ is constant, $\mu\equiv0$, and~(\ref{eq:explicit}) returns
$f = \cos$, again up to a quadratic.
\end{example}

\begin{example}
[Powers of a norm]\label{ex:normpow}Let $X$ be a real normed space with norm
$\|\cdot\|$, take $G = X$, and let $f(x) = \|x\|^{s}$ with $0 < s < 2$. As $f$ is
even, $g = \mathsf{Q}[f]$ has $g(0, x) = 0$ and $\mu\equiv0$, and on the diagonal
$g(x, x) = \|2 x\|^{s} - 4\|x\|^{s} = (2^{s} - 4)\|x\|^{s}$. The datum is
homogeneous of degree $s$, so $\|g(2^{n-1}x, 2^{n-1}y)\| = 2^{(n-1)s}\|g(x,
y)\|$, and~(\ref{eq:conv}) holds precisely because $s < 2$.
Formula~(\ref{eq:explicit}) then gives
\[
f(x) = -(2^{s} - 4)\|x\|^{s}\sum_{n\geq1} 4^{-n} 2^{(n-1)s} = \|x\|^{s},
\]
recovering $f$. For $s = 1$ this is the norm itself, of which
Example~\ref{ex:abs} is the one-dimensional case; in higher dimensions the
off-diagonal values are nontrivial, so the propagation does substantive work.
\end{example}

\begin{example}
[The parallelogram-law defect]\label{ex:defect}Let $\|\cdot\|$ be a norm not
induced by an inner product---an $\ell^{p}$ norm with $p\neq2$, say---and take
$f(x) = \|x\|^{2}$. By the Jordan--von Neumann theorem the defect $g =
\mathsf{Q}[f]$ is not identically zero, and two things coincide. First, $g$ is
homogeneous of degree $2$, the critical exponent:
\[
\sum_{n\geq1} 4^{-n}\bigl\|g(2^{n-1}x, 2^{n-1}y)\bigr\| =
\sum_{n\geq1}\tfrac{1}{4}\|g(x, y)\|,
\]
so~(\ref{eq:conv}) fails. Second, by homogeneity the diagonal vanishes, $g(x, x)
= \|2 x\|^{2} - 4\|x\|^{2} = 0$, so the diagonal-driven
formula~(\ref{eq:explicit}) would return $f\equiv0$ although $g\not\equiv0$. Thus
$g$ is a quasi-symmetric cocycle in $\Image\mathsf{Q}$, with $\|\cdot\|^{2}$ a
solution unique up to $\Ker\mathsf{Q}$, lying on the exact boundary where the
construction fails; the vanishing of the diagonal shows why~(\ref{eq:conv}) is
indispensable, the off-diagonal values being lost to the diagonal once the
summability that propagates them is gone. This sits at the critical exponent,
sharper than the super-quadratic datum discussed in
Section~\ref{sec:cohomology}.
\end{example}

\begin{example}
[A discrete group]\label{ex:Z}Take $G = \mathbb{Z}$ and $X = \mathbb{R}$, with
the datum
\[
g(x, y) = |x + y| + |x - y| - 2 |x| - 2 |y|,
\]
the restriction to the integer lattice of the datum of
Example~\ref{ex:abs}. It is a quasi-symmetric cocycle ($g = \mathsf{Q}[h]$ with
$h = |\cdot|$), and its linear growth makes~(\ref{eq:conv}) hold. Here $g(0,
m) = 0$, so $\mu\equiv0$, and since $g_{e}(t, t) = -2 |t|$,
formula~(\ref{eq:explicit}) gives
\[
f(m) = -\sum_{n\geq1} 4^{-n}\bigl(-2\,|2^{n-1}m|\bigr) = |m|\sum_{n\geq1}
2^{-n} = |m|,
\]
recovering $h$. Two features distinguish this from the real-variable examples.
First, no topology is available, so the regularity hypotheses of the smooth
theory are vacuous and the only input is the summability~(\ref{eq:conv}).
Second, $\mathbb{Z}$ is not $2$-divisible---the odd integers are not
doubles---yet the construction succeeds, because both its stages use only the forward
orbit $m, 2 m, 4 m, \dots$ and never halve.
The same computation applies verbatim on $\mathbb{Z}^{d}$ or on any finitely
generated free abelian group; a linear term added to $h$ is recovered through
the odd component $\mu$.
\end{example}

\begin{example}
[Finite abelian groups]\label{ex:finite}Let $G$ be any finite abelian group, $X
= \mathbb{R}$, and $g = \mathsf{Q}[h]$ for an arbitrary $h\colon G\to\mathbb{R}$.
Three features stand out. Condition~(\ref{eq:conv}) holds automatically, $g$
being bounded, so no growth or regularity hypothesis is needed and membership in
$\Image\mathsf{Q}$ becomes purely algebraic. The solution is unique: a quadratic
$f(x) = B(x, x)$ with $B$ biadditive into the torsion-free $\mathbb{R}$ has
$N\, B(x, y) = B(N x, y) = 0$ on a group of order $N$, since $N x = 0$ there,
forcing $B\equiv0$, so $\Ker\mathsf{Q} = \{0\}$ and~(\ref{eq:explicit}) returns
$h$ itself. Finally the orbit $2^{n}x$ is eventually periodic, so the series
in~(\ref{eq:explicit}) is a finite geometric sum: on $\mathbb{Z}/N\mathbb{Z}$
with $N$ odd, multiplication by $2$ is invertible and the orbits are full
cycles, whereas on $\mathbb{Z}/2^{k}\mathbb{Z}$ it is nilpotent and the orbit
collapses to $0$, the reconstruction holding either way. This is the algebraic
counterpart of the non-quasi-symmetric cocycle exhibited on
$(\mathbb{Z}/N\mathbb{Z})^{2}$ in Section~\ref{sec:cohomology}.
\end{example}

\begin{remark}
\label{rem:compact}More generally, on any compact abelian group with $X =
\mathbb{R}$ and $g$ continuous, $g$ is bounded and~(\ref{eq:conv}) holds for
free, finite groups being the discrete case. The $2$-adic integers
$\mathbb{Z}_{2}$ are striking: $\|2^{n}x\|_{2}\to0$, so the doubling orbit tends
to $0$, the exact reverse of the real line, where $2^{n}x\to\infty$. The
defining series of~(\ref{eq:explicit}) then converges by mere boundedness of its
terms, the construction propagating toward $0$ rather than toward infinity.
\end{remark}

\begin{example}
[An infinite-dimensional target]\label{ex:banach}Take $G = \mathbb{Z}$ and $X =
\ell^{2}(\mathbb{N})$ with standard basis $(e_{k})$, and let
\[
f(m) = \sum_{k\geq1} 2^{-k}\cos(k m)\, e_{k}.
\]
Then $f(m)\in\ell^{2}$, with $\|f(m)\|\leq(\sum_{k\geq1} 4^{-k})^{1/2}$, and the
datum $g = \mathsf{Q}[f]$ is bounded, so~(\ref{eq:conv}) holds trivially. The
axial trace $g(0, m) = -2\sum_{k\geq1} 2^{-k}e_{k}$ is constant, hence
$\mu\equiv0$, and formula~(\ref{eq:explicit}) reconstructs $f$, its partial
sums converging in $\ell^{2}$ precisely because $X$ is complete. This is where
completeness of $X$ enters: it turns the formal series of the construction into a
convergent one. Nothing is special to $\ell^{2}$; any Banach space and
any datum satisfying~(\ref{eq:conv}) would serve. For instance the same $f$ lies
in $c_{0}$, which carries no inner product and is not reflexive, so neither
property is used.
\end{example}

\section{Application to Ulam stability}

\label{sec:stability}

\emph{Hyers--Ulam stability} asks whether an approximate solution of a functional
equation must lie close to an exact one. The question was raised by
Ulam~\cite{Ulam1960} in the context of metric groups and answered by
Hyers~\cite{Hyers1941} for the additive Cauchy equation in Banach spaces: if
$\|f(x + y) - f(x) - f(y)\|\leq\varepsilon$ for all $x, y$, then there is a
unique additive map $A$ with $\|f - A\|\leq\varepsilon$. The same scheme applies
to the quadratic equation~(\ref{eq:JvN}), where one asks whether
$\|\mathsf{Q}[f]\|\leq\varepsilon$ forces $f$ to be near a quadratic function;
this was established by Skof~\cite{Skof1983} and Cholewa~\cite{Cholewa1984}.

Two relaxations of the uniform bound are by now standard.
Rassias~\cite{Rassias1978} allowed power growth $\theta(\|x\|^{p} + \|y\|^{p})$
with $p$ below a critical exponent, so that the admissible error may grow with
the variables---for the quadratic equation this is due to
Czerwik~\cite{Czerwik1992}. G\u{a}vru\c{t}a~\cite{Gavruta1994}, building on the iterative
\emph{direct method} of Forti~\cite{Forti1980,Forti1987}, then replaced power
growth by an \emph{arbitrary} control function $\varphi$ subject only to a
summability condition. We work at this last level of generality; see the
monograph~\cite{BrzdekPopaRasaXu2018} for a systematic account.

The construction of Section~\ref{sec:descent} is precisely this direct method:
applied to an \emph{approximate} solution it yields the stability of the quadratic
equation in the Forti--G\u{a}vru\c{t}a sense. Throughout this section $\varphi\colon
G\times G\to[0, \infty)$ satisfies the \emph{Forti--G\u{a}vru\c{t}a \mbox{condition}}
\begin{equation}
\label{eq:gavruta}\widetilde{\varphi}(x, y) := \sum_{n\geq1} 4^{-n}\varphi
(2^{n-1}x, 2^{n-1}y) < \infty,\qquad x, y\in G;
\end{equation}
this is exactly the two-variable condition~(\ref{eq:conv}) of the inverse
problem, with $\varphi$ in place of $\|g(\cdot, \cdot)\|$.

\begin{theorem}
\label{thm:stab}If $f\colon G\to X$ satisfies $\|\mathsf{Q}[f](x, y)\|\leq
\varphi(x, y)$ for all $x, y\in G$, then there is a unique quadratic $q\colon
G\to X$ with
\begin{equation}
\label{eq:stabound}\|f(x) - q(x)\|\leq\widetilde{\varphi}(x, x) + \tfrac{1}{6}
\varphi(0, 0),\qquad x\in G,
\end{equation}
namely $q(x) = \lim\limits_{n\to\infty} 4^{-n} f(2^{n}x)$, $x\in G$.
\end{theorem}

\begin{proof}
Set $\tilde{f} := f - \tfrac{1}{3} f(0)$; the shift is chosen so that
$\tilde{f}(2 x) - 4 \tilde{f}(x) = \mathsf{Q}[f](x, x)$ for all $x$, the constant
absorbing the term $f(0)$ in $\mathsf{Q}[f](x, x) = f(2 x) + f(0) - 4 f(x)$. Thus
$\|\tilde{f}(2 x) - 4 \tilde{f}(x)\|\leq\varphi(x, x)$, and
\begin{multline*}
\bigl\|4^{-(n+1)}\tilde{f}(2^{n+1}x) - 4^{-n}\tilde{f}(2^{n}x)\bigr\|\\
= 4^{-(n+1)}\bigl\|\tilde{f}(2^{n+1}x) - 4 \tilde{f}(2^{n}x)\bigr\|\leq
4^{-(n+1)}\varphi(2^{n}x, 2^{n}x)
\end{multline*}
for all $n\geq0$.
The right-hand side is summable by~(\ref{eq:gavruta}), so
$(4^{-n}\tilde{f}(2^{n}x))_{n}$ is Cauchy and, $X$ being complete, converges to a
limit $q(x) = \lim\limits_{n\to\infty} 4^{-n}\tilde{f}(2^{n}x) = \lim\limits_{n\to\infty} 4^{-n} f(2^{n}x)$, the
constant $\tfrac{1}{3} f(0)$ vanishing in the limit. Summing the telescoping
bound over $n\geq0$ gives $\|\tilde{f}(x) - q(x)\|\leq\sum_{n\geq1}
4^{-n}\varphi(2^{n-1}x, 2^{n-1}x) = \widetilde{\varphi}(x, x)$, so
\[
\|f(x) - q(x)\|\leq\widetilde{\varphi}(x, x) + \tfrac{1}{3}\|f(0)\|.
\]
Since $\|\mathsf{Q}[f](0, 0)\| = \|{-2} f(0)\|\leq\varphi(0, 0)$ forces
$\|f(0)\|\leq\tfrac{1}{2}\varphi(0, 0)$, this is~(\ref{eq:stabound}).\smallskip

\noindent\emph{$q$ is quadratic.} For fixed $x, y$,
\[
\|\mathsf{Q}[\tilde{f}](2^{n}x, 2^{n}y)\| = \|\mathsf{Q}[f](2^{n}x, 2^{n}y) +
\tfrac{2}{3} f(0)\|\leq\varphi(2^{n}x, 2^{n}y) + \tfrac{2}{3}\|f(0)\|,
\]
so $4^{-n}\mathsf{Q}[\tilde{f}](2^{n}x, 2^{n}y)\to0$ by~(\ref{eq:gavruta});
hence $\mathsf{Q}[q](x, y) = \lim\limits_{n\to\infty}
4^{-n}\mathsf{Q}[\tilde{f}](2^{n}x, 2^{n}y) = 0$, so $q\in\Ker\mathsf{Q}$.\smallskip

\noindent\emph{Uniqueness.} If $q, q_{1}$ are quadratic and both
satisfy~(\ref{eq:stabound}), then $r := q - q_{1}$ is quadratic, so $r(2^{n}x) =
4^{n} r(x)$ and
\[
\|r(x)\| = 4^{-n}\|r(2^{n}x)\|\leq 4^{-n}\bigl(2\widetilde{\varphi}(2^{n}x,
2^{n}x) + \tfrac{1}{3}\varphi(0, 0)\bigr).
\]
As $4^{-n}\widetilde{\varphi}(2^{n}x, 2^{n}x)$ is the tail of the convergent
series~(\ref{eq:gavruta}) and tends to $0$, we get $r(x) = 0$.
\end{proof}

The two-variable hypothesis~(\ref{eq:gavruta}) is needed only to make $q$
quadratic; the bound~(\ref{eq:stabound}) itself involves $\varphi$ on the
diagonal alone. When $\varphi$ is dominated by its diagonal the distinction
disappears.

\begin{corollary}
\label{cor:diagdom}Suppose $\varphi(x, y)\leq\tfrac{1}{2}\bigl(\varphi(x, x) +
\varphi(y, y)\bigr)$ for all $x, y$, and that $\sum_{n\geq1} 4^{-n}\varphi
(2^{n-1}x, 2^{n-1}x) < \infty$ for all $x$. Then the conclusion of
Theorem~\textnormal{\ref{thm:stab}} holds. This covers the classical controls:
a constant $\varepsilon$; the Rassias control~\cite{Rassias1978} $\theta(\|x\|
^{p} + \|y\|^{p})$ with $p < 2$; and the product control $\theta\|x\|^{p}\|y\|
^{p}$ with $p < 1$, on a normed space $G$.
\end{corollary}

\begin{proof}
The domination gives $4^{-n}\varphi(2^{n}x, 2^{n}y)\leq\tfrac{1}{2}\bigl(4^{-n}
\varphi(2^{n}x, 2^{n}x) + 4^{-n}\varphi(2^{n}y, 2^{n}y)\bigr)$, and each term on
the right tends to $0$ since the diagonal series converges; thus the only use
of~(\ref{eq:gavruta}) in the proof of Theorem~\ref{thm:stab} is supplied by the
diagonal hypothesis. The three controls satisfy the domination: the constant and Rassias controls
with equality, and the product control by the inequality of arithmetic and
geometric means.
\end{proof}

\begin{remark}
\label{rem:classical}For the constant control $\varphi\equiv\varepsilon$ above,
the bound~(\ref{eq:stabound}) becomes $\|f(x) - q(x)\|\leq\tfrac{1}{2}
\varepsilon$, since $\widetilde{\varphi}(x, x) = \sum_{n\geq1} 4^{-n}\varepsilon
= \tfrac{1}{3}\varepsilon$ and $\tfrac{1}{6}\varphi(0, 0) =
\tfrac{1}{6}\varepsilon$. This recovers the classical Hyers--Ulam stability of
the quadratic equation, due to Skof~\cite{Skof1983} and
Cholewa~\cite{Cholewa1984}, with the optimal Ulam constant $\tfrac{1}{2}$.
\end{remark}

The two halves of the paper meet here: the stability quadratic is exactly the
quadratic defect of $f$ relative to the inverse problem.

\begin{remark}
\label{rem:link}With $f$, $\varphi$, $q$ as in
Theorem~\textnormal{\ref{thm:stab}}, the function $f_{1} := f - q$ satisfies
$\mathsf{Q}[f_{1}] = \mathsf{Q}[f]$, since $q\in\Ker\mathsf{Q}$. Thus $f = q +
f_{1}$ decomposes $f$ into its quadratic part $q$ and a solution $f_{1}$ of the
inhomogeneous equation~(\ref{eq:inhom}) for the datum $g = \mathsf{Q}[f]$, which
satisfies~(\ref{eq:cocycle}) and~(\ref{eq:qsym}) by Theorem~\ref{thm:necessity}
and, being bounded by $\varphi$, the hypothesis~(\ref{eq:conv}). Hence
Theorem~\ref{thm:inverse} applies, and $f_{1}$ differs from the
solution~(\ref{eq:explicit}) by a quadratic.
\end{remark}

\section{A cohomological interpretation}

\label{sec:cohomology}

The cocycle identity admits a cohomological interpretation that places the linear and
quadratic inverse problems in parallel. For the Cauchy operator $\mathsf{L}$,
the coboundary map $f\mapsto\mathsf{L}[f]$ and the cocycle defect fit into the
standard bar complex; symmetric $2$-cocycles modulo $2$-coboundaries form
$\mathrm{Ext}^{1}_{\mathbb{Z}}(G, X)$, which vanishes because $X$, a real vector
space, is divisible and hence an injective $\mathbb{Z}$-module
(see~\cite[Sec.~20--24]{Fuchs1970}). This is the
structural reason behind the unconditional characterization of $\Image
\mathsf{L}$ recalled in the introduction: \emph{no} regularity or growth
hypothesis is required.

The operator $\mathsf{Q}$ admits the analogous pairing $\mathsf{D}\circ\mathsf{Q}
= 0$ (Theorem~\ref{thm:necessity}), so $\Image\mathsf{Q}\subseteq\Ker\mathsf{D}$,
and the quotient $\Ker\mathsf{D}\,/\,\Image\mathsf{Q}$ measures the failure of
surjectivity onto cocycles. Two features distinguish it from the linear case.
First, $\Ker\mathsf{D}$ must be cut down by quasi-symmetry: the difference $g(x, y) - g(y, x)$ is an independent
invariant, vanishing modulo the
axial trace exactly on the quasi-symmetric cocycles. Second---and unlike the
divisible vanishing available for $\mathsf{L}$---we know of no purely algebraic
reason why a quasi-symmetric cocycle should be a coboundary. What
Theorems~\ref{thm:necessity} and~\ref{thm:inverse} establish is that the
Forti--G\u{a}vru\c{t}a condition~(\ref{eq:conv}) suffices to annihilate the
obstruction: through the doubling relation~(\ref{eq:star}), every
quasi-symmetric cocycle satisfying~(\ref{eq:conv}) is a coboundary. The
condition is sufficient but not necessary, and it can fail already at quadratic
growth: the parallelogram-law defect of a non-Hilbertian norm
(Example~\ref{ex:defect}) lies in $\Image\mathsf{Q}$ yet
violates~(\ref{eq:conv}), as does the coboundary of any super-quadratic
function, say $f(x) = |x|^{3}$ on $\mathbb{R}$; in each case the defining series
of the inverse diverges. So the construction reaches precisely those
quasi-symmetric cocycles that satisfy~(\ref{eq:conv}). Whether
every quasi-symmetric cocycle, irrespective of growth, is a
coboundary---equivalently, whether $\Ker\mathsf{D}$ restricted to
quasi-symmetric cocycles coincides with $\Image\mathsf{Q}$ over an arbitrary
abelian group---we leave as an open problem.

A second point concerns quasi-symmetry itself, and here the situation is
settled: it is genuinely independent of the cocycle identity. On $G =
\mathbb{Z}^{2}$ the function
\[
g\bigl((x_{1}, x_{2}), (y_{1}, y_{2})\bigr) = x_{2}^{2}\,y_{1}\,y_{2}
\]
satisfies the cocycle identity~(\ref{eq:cocycle}) identically---a direct
substitution gives $\mathsf{D}[g] = 0$---yet it is not quasi-symmetric: since
$g$ vanishes whenever its first argument is $0$, the right-hand side
of~(\ref{eq:qsym}) is identically zero, while the left-hand side $g(x, y) - g(y,
x) = x_{2}^{2}y_{1}y_{2} - y_{2}^{2}x_{1}x_{2}$ equals $1$ at $x = (0, 1)$, $y =
(1, 1)$. Hence quasi-symmetry cannot be dropped from
Theorem~\ref{thm:necessity}; it must be imposed alongside~(\ref{eq:cocycle}),
exactly as symmetry must be imposed alongside the additive cocycle in the linear
case by Erd\H{o}s~\cite{Erdos1959}. The construction uses two independent
generators of $\mathbb{Z}^{2}$; the same $g$, read modulo $N$, is a
non-quasi-symmetric cocycle on $(\mathbb{Z}/N\mathbb{Z})^{2}$ as well, so the
phenomenon is one of rank rather than of infinitude.

Independence is the generic situation, not a universal one: whenever the
homogeneous equation has sufficient uniqueness, the cocycle identity does force
symmetry. The argument uses only the symmetry of $T[\,\cdot\,]$ and a uniqueness
property, so it applies to any operator with symmetric values; for the additive
cocycle it is the one given by Acz\'{e}l and presented by Erd\H{o}s
in~\cite{Erdos1959}.

\begin{proposition}
[A symmetry criterion]\label{prop:symmetry-criterion}Let $V_{1}\subseteq
\mathcal{F}(G, X)$ be a linear subspace and put
\[
V_{2} = \bigl\{\, g\in\mathcal{F}(G^{2}, X) : g(t, \cdot)\in V_{1}\text{ and }
g(\cdot, t)\in V_{1}\text{ for all }t\in G\,\bigr\}.
\]
Let $T\colon V_{1}\to V_{2}$ be linear with $T[f]$ symmetric for all $f\in
V_{1}$, and suppose that
\begin{equation}
\label{eq:uniqueness}\text{for all }x\in G\setminus\{0\}\text{, the only }
f\in V_{1}\text{ with }T[f] = 0\text{ and }f(x) = 0\text{ is }f = 0.
\end{equation}
Then every $g\in V_{2}$ satisfying
\begin{equation}
\label{eq:itercocycle}T[g(x, \cdot)](y, z) = T[g(\cdot, z)](x, y),\qquad x, y,
z\in G,
\end{equation}
is symmetric.
\end{proposition}

\begin{proof}
Fix $x\in G\setminus\{0\}$ and set $f := g(x, \cdot) - g(\cdot, x)$, which lies
in $V_{1}$ because $g\in V_{2}$. Writing $C$ for the
hypothesis~(\ref{eq:itercocycle}) and $S$ for the symmetry of $T[\,\cdot\,]$, for
all $y, z\in G$,
\begin{align*}
T[g(x, \cdot)](y, z)  &  \overset{C}{=} T[g(\cdot, z)](x, y) \overset{S}{=}
T[g(\cdot, z)](y, x) \overset{C}{=} T[g(y, \cdot)](x, z)\\
&  \overset{S}{=} T[g(y, \cdot)](z, x) \overset{C}{=} T[g(\cdot, x)](y, z).
\end{align*}
Hence $T[f](y, z) = 0$ for all $y, z$, so $f\in\Ker T$. Since $f(x) = g(x, x) -
g(x, x) = 0$, condition~(\ref{eq:uniqueness}) gives $f = 0$, that is, $g(x, t) =
g(t, x)$ for all $t\in G$. This holds for all $x\ne0$, so $g(x, y) = g(y, x)$
whenever $x\ne0$ and---exchanging the two arguments---also whenever $y\ne0$; the
case $x = y = 0$ is trivial. Hence $g$ is symmetric.
\end{proof}

Only the symmetry of $T[\,\cdot\,]$ and the uniqueness~(\ref{eq:uniqueness}) are
used; the cocycle identity for $T$ itself (the iterated identity of
Remark~\ref{rem:operatorial}) plays no role. The hypotheses hold for both
operators of this paper. For $\mathsf{L}$ every value $\mathsf{L}[f](x, y)$ is
symmetric, as~(\ref{eq:Ldef}) is symmetric in $x, y$, so the criterion applies
with $V_{1} = \mathcal{F}(G, X)$; for $G = X = \mathbb{R}$ and $V_{1}$ the
continuous
functions, a continuous solution of $\mathsf{L}[f] = 0$ is additive, hence
linear, and vanishing at one nonzero point forces it to vanish,
so~(\ref{eq:uniqueness}) holds. The proposition then states that every
continuous solution of the additive cocycle on $\mathbb{R}$ is
symmetric---equivalently, by the divisible vanishing above, a
coboundary---which is the result of Erd\H{o}s~\cite{Erdos1959}.

For $\mathsf{Q}$, identity~(\ref{eq:el-swap}) shows that $\mathsf{Q}[f]$ is
symmetric exactly when $f$ is even; take $V_{1}$ to be the even functions, so
that $V_{2}$ consists of the functions even in each variable. For $G = X =
\mathbb{R}$ and $V_{1}$ the continuous even functions,~(\ref{eq:uniqueness})
holds: a continuous even solution of $\mathsf{Q}[f] = 0$ is a quadratic $f(x) =
c\,x^{2}$, which vanishes at a nonzero point only if $c = 0$. Given a continuous
$g$ satisfying~(\ref{eq:cocycle}) on $\mathbb{R}$, the even reduction $g_{e}$
of~(\ref{eq:gedef}) is continuous, even in each variable, and again
satisfies~(\ref{eq:cocycle}) by Proposition~\ref{prop:residual}; the criterion
applies to $g_{e}$, so $g_{e}$ is symmetric and $g$ is quasi-symmetric by
Proposition~\ref{prop:qsym-equiv}. Thus on $\mathbb{R}$ every continuous cocycle
is quasi-symmetric. Continuity may here be relaxed to any of the mild regularity
hypotheses recalled in the introduction---Lebesgue measurability, or boundedness
on a set of positive measure---each still forcing the relevant kernel element to
be linear, respectively a quadratic monomial; this is the weakening already
signalled by Erd\H{o}s in~\cite{Erdos1959}. The independence exhibited on
$\mathbb{Z}^{2}$ accordingly disappears under any such regularity, not under
continuity alone.

\section*{Statements and Declarations}

\subsection*{Funding}
No funding was received to assist with the preparation of this manuscript.

\subsection*{Competing Interests}
The authors have no competing interests to declare that are relevant to the
content of this article.

\subsection*{Use of AI Tools}
The authors used the large language model Claude (Anthropic) to review this
manuscript: to check its language, wording and structure, to look for errors in
its arguments, and to verify its references and the attribution of results. The
results and their proofs are the authors' own, and the authors take full
responsibility for the content of the article.

\bibliographystyle{spmpsci}
\bibliography{inverse-quadratic-groups}

\end{document}